\documentclass[a4paper,pdftex,reqno,10pt]{amsart}
\usepackage{hyperref}
\usepackage{amsmath,amssymb}
\usepackage{url}
\usepackage[textwidth=160mm, textheight=245mm]{geometry}

\begin{document}
\hyphenation{equiva-lence}

\date{\today}

\title{On minimum sum representations for weighted voting games}
\author{Sascha Kurz}
\address{Department of Mathematics, Physics, and Computer Science, University of Bayreuth, 95440 Bayreuth, Germany\\Tel.: +49-921-557353, Fax: +49-921-557352}
\email{sascha.kurz@uni-bayreuth.de}
\begin{abstract}
  A proposal in a weighted voting game is accepted if the sum of the (non-negative) weights of the ``yea'' voters is at
  least as large as a given quota. Several authors have considered representations of weighted voting games with minimum
  sum, where the weights and the quota are restricted to be integers. In \cite{root} the authors have classified all
  weighted voting games without a unique minimum sum representation for up to $8$ voters. 
  
  Here we exhaustively classify all weighted voting games consisting of $9$~voters which do not admit a unique
  minimum sum integer weight representation.

  \bigskip

\end{abstract}

%
%
%
%
%
%
%

\maketitle

\section{Introduction}

\noindent
Consider a yes-no voting system for a set of~$n$ voters. The acceptance of a proposal should depend on the subset of ``yea'' voters. 
There are different concepts for the set of required features of a voting system. One is that of a weighted voting game\footnote{Other aliases are weighted (majority) games or threshold functions.}. Here we are given non-negative voting weights $w_i\in\mathbb{R}_{\ge 0}$ for the voters and a quota $q\in\mathbb{R}_{>0}$. A proposal is accepted iff $\sum\limits_{i\in Y}w_i\ge q$, where $Y$ is the set of voters which are in favor of the proposal. 
In \cite{root} the authors restrict the weights and the quota to be integers
and ask for weight representations with minimum sum. They have shown that for at most $7$~voters these representations are unique.
For $8$ voters there are exactly $154$~weighted majority games with two minimum sum integer weight representations. If
one requires that equally desirable voters obtain equal weights one speaks of a normalized representation. Indeed all
these $154$~weighted games admit a unique normalized minimum sum integer weight representation. They also give some examples
consisting of $10$ or more voters which do not admit such a unique normalized representation. Recently the same authors
have presented such examples for $9$~voters in \cite{root2}. We would like to remark that minimum integer weight representations of weighted voting
games are also used as a solution concept for cooperative transferable utility games, see \cite{0871.90128}.

\subsection{Related results}
\noindent
Isbell, see \cite{isbell}, found an example of $12$~(unsymmetrical) voters without a unique minimum sum representation. Smaller examples for $9$, $10$, and $11$~voters are given 
in \cite{root}, \cite{root2}. The enumeration of weighted voting games dates back to at least 1962 \cite{0105.12002}, where up to $6$~voters are treated. For $n=7,8$ voters we refer 
e.g.\ to \cite{0233.94016}, \cite{0205.17805}, \cite{0841.90134}. Bart de Keijzer presents a gra\-ded poset for weighted voting games in his master thesis \cite{keijzer}, see also 
\cite{keijzer2}.\footnote{We remark that the counts for weighted voting games with $6\le n\le 8$ voters are wrongly stated in \cite{keijzer}. The causative buglet is fixed by now 
(personal communication).} To our knowledge, the number of weighted voting games for $9$~voters has not been published before.\footnote{We would like to mention the (unpublished) 
diploma thesis \cite{da_tautenhahn} containing the enumeration for $9$~voters (without using the results from this paper; compare footnote \ref{fn_da_tautenhahn}).}

\subsection{Our contribution}
\noindent
We exhaustively classify all weighted majority games consisting of $9$~voters which do not admit a unique (normalized) minimum sum integer weight representation. Within this context 
we also enumerate the number of weighted voting games for $9$~voters. The new enumeration results are possible due to the newly introduced concept of partial complete simple games 
and efficient implementations.

\subsection{Outline of the Paper}
\noindent
In Section~\ref{sec_subclass} we very briefly\footnote{For a more extensive introduction we refer to \cite{0943.91005}.} state the basic definitions and facts of weighted voting games 
in the context of simple games. In Section~\ref{sec_partial} we describe an approach how to exhaustively enumerate weighted voting games without generating the whole set of complete 
simple games. The minimum sum representations are treated in Section~\ref{sec_min_rep}.

\section{Weighted voting games as a subclass of (complete) simple games}
\label{sec_subclass}
\noindent
Let $N=\{1,\dots,n\}$ be a set of $n$ voters. An example for a weighted voting game with $n=4$ voters, as defined in the introduction, is given by quota $q=3$ and the weights $w_1=2$, $w_2=w_3=1$, $w_4=0$. By comparing $\sum_{i\in U} w_i$ with $q$ we can decide for each coalition $U\subseteq N$ whether it is winning (the proposal gets accepted) or losing. So being more general a pair $(N,\chi)$ is called a simple game if $\chi$ is a characteristic function of the subsets of $N$ with $\chi(\emptyset)=0$, $\chi(N)=1$, and $\chi(U')\le\chi(U)$ for all $U'\subseteq U$. For the broad variety of applications of simple games we quote Taylor and Zwicker \cite{0943.91005}: ``Few structures arise in more contexts and lend themselves to more diverse interpretations than do simple games.'' The subset-minimal winning coalitions of our example are given by $\{1,2\}$, $\{1,3\}$ and uniquely characterize the simple game.

A well studied subclass (and superclass of weighted voting games) arises from Isbell's desirability relation \cite{0083.14301}: We write $i\sqsupset j$ for two voters $i,j\in N$ iff we have $\chi\Big(\{i\}\cup U\backslash\{j\}\Big)\ge\chi (U)$ for all $j\in U\subseteq N\backslash\{i\}$. A pair $(N,\chi)$ is called complete simple game if it is a simple game and the binary relation $\sqsupset$ is a total preorder. We abbreviate $i\sqsupset j$, $j\sqsupset i$ by $i~\square~j$. In our example we have $1\sqsupset2~\square~3\sqsupset 4$ ($2\not\sqsupset 1$ and $4\not\sqsupset 3$).

In the following we restrict our considerations onto symmetry classes of simple games, i.e.\ we consider games arising from renumbering the set of players as equivalent, see \cite{complete_simple_games} for more details, and assume $1\sqsupset 2 \sqsupset \dots \sqsupset n$ (corresponding to $w_1\ge\dots\ge w_n$ for weighted voting games). We write coalitions $U\subseteq N$ as characteristic vectors\footnote{Another representation, which takes the possible symmetry of voters into account, is described e.g.\ in \cite{complete_simple_games}. There are several applications where this representation is more convenient, see e.g.\ \cite{dedekind}.} $u\in\{0,1\}^n$. The two subset-minimal winning coalitions of our example correspond to $(1,1,0,0)$ and $(1,0,1,0)$. We call $u$ a winning coalition iff $\chi(U)=1$, otherwise we call $u$ a losing coalition. To have a compact representation for complete simple games we need another partial ordering: For two coalitions $u=\left(u_1,\dots,u_n\right), v=\left(v_1,\dots,v_n\right)\in\{0,1\}^n$ we write $u\preceq v$ iff we have $\sum\limits_{i=1}^{k} u_i \le \sum\limits_{i=1}^{k} v_i$ for all $1\le k\le n$. We abbreviate $u\preceq v$, $u\neq v$ by $u\prec v$ and $u \npreceq v$, $u \nsucceq v$ by $u\bowtie v$. An example is given by $(1,0,1,0)\prec(1,1,0,0)$. For a complete simple game we denote by $W$ all, so-called shift-minimal, winning coalitions which are minimal with respect to $\preceq$. Similarly we denote by $L$ the set of all, so-called shift-maximal, losing coalitions which are maximal with respect to $\preceq$. Due to $\chi(u)\le \chi(v)$ for all $u\preceq v$ each complete simple game is uniquely characterized by either $W$ or $L$, see e.g.\ \cite{complete_simple_games}. In our example the set of shift-minimal winning coalitions is given by $W=\{(1,0,1,0)\}$ implying that $(1,1,0,0)\succ(1,0,1,0)$ is also a winning coalition.

Let $G_n$ be the graph consisting of vertex set $\{0,1\}^n$ and edges $\{u,v\}$ for all $u,v\in\{0,1\}^n$ with $u\bowtie v$. The non-empty cliques of $G_n$ are in bijection to the sets $W$ of minimal winning coalitions and thereby to complete simple games for $n$ voters \cite{dedekind}. Applying the software package \texttt{cliquer} \cite{cliquer}, \cite{1019.05054}, using a customary personal computer, we have the counts in Table~\ref{table_csg_cliquer} (coinciding with those from \cite{root2}).

There are several approaches how to check whether a complete simple game is weighted, for an over\-view see e.g.\ \cite{root}, \cite{0943.91005}. Minimum sum representations of weighted voting games are e.g.\ in bijection to the solutions of the following integer linear program
\begin{align}
  \min && \sum_{i=1}^n w_i\\
  s.t.{} &&  \sum_{i=1}^n u_iw_i \ge 1+\sum_{i=1}^n v_iw_i &&\forall u\in W, v\in L,\\
  && w_i\ge w_{i-1} && \forall 1\le i\le n-1,\\
  && w_i \in\mathbb{Z}_{\ge 0} && \forall 1\le i\le n,
\end{align}
which indeed was used in \cite{root}. A complete simple game is weighted iff a solution of this integer linear program or its linear relaxation exists. As mentioned in \cite{root}, generating all $284\cdot 10^9$ complete simple games for $9$~voters and afterwards solving the corresponding (integer) linear programs would be too time-consuming.

\begin{table}[htp]
\begin{center}
\begin{tabular}{crrrrrrrrr}
\hline
$\mathbf{n}$&$1$&$2$&$3$&$4$&$5$&$6$&
$7$&$8$&$9$\\
$\mathbf{\#}$&$1$&$3$&$8$&$25$&$117$&$1171$&
$44313$&$16175188$&$284432730174$\\
$\textbf{time}$&&&&&&&&$1~s$&$44~m$\\
\hline
\end{tabular}
\caption{Complete simple games for $n$~voters up to symmetry ($1\sqsupset 2 \sqsupset \dots \sqsupset n$).}
\label{table_csg_cliquer}
\end{center}
\end{table}

\section{Partial complete simple games and similar linear programs}
\label{sec_partial}
\noindent
Suppose we want to exhaustively generate the complete simple games for a given number of voters using an orderly generation approach, see \cite{0392.05001}. To be more precise, we start with an empty set $W$ and add shift-minimal winning coalitions, which are decreasing with respect to the lexicographical ordering $\le_{\text{lex}}$\footnote{Here one may also read the coalitions as integers written in their binary expansion and use the ordinary ordering $\le$ of integers.}, step by step. In the corresponding generation tree for a node $W$ all successors $W'$ fulfill $W\subseteq W'$ and $y<_{\text{lex}}x$ for all $x\in W$, $y\in W'\backslash W$. Now let $L$ and $L'$ be the sets of shift-maximal losing coalitions corresponding to $W$ and $W'$. For a given set $W\neq\emptyset$ with lexicographically smallest element $w$ we define
\begin{eqnarray*}
  C:=W\cup \left\{x\in\{0,1\}^n\mid x<_{\text{lex}}w,\nexists y\in W:x\prec y \right\}
\end{eqnarray*}
and
\begin{eqnarray*}
  \tilde{L}&:=&\Big\{v\in\{0,1\}^n\mid\exists u\in W:\,v\prec u\Big\}\,
  \cup\,\Big\{v\in\{0,1\}^n\mid\nexists u\in C
  :\,v\succeq u\Big\}
\end{eqnarray*}
and condense the maximal coalitions of $\tilde{L}$ to a set $\hat{L}\subseteq\tilde{L}$. For $v'\in\Big\{v\in\{0,1\}^n\mid\exists u\in W:\,v\prec u\Big\}$ choose $u'\in W$ with $v'\prec u'$. Since $u'$ is a shift-minimal winning coalition in the complete simple game characterized by $W'$ the coalition $u'$ must be a losing coalition. Using the orderly generation approach and the fact that $W'$ is an antichain we conclude $W'\subseteq C$. A coalition $v$ can only be winning if there is an $u\in C\supseteq W'$ with $v\succeq u$. Thus we have $\hat{L}\subseteq\tilde{L}\subseteq L$ and $\hat{L}\subseteq\tilde{L}\subseteq L'$ for all successors of $W$.

Now we may use the partial sets $W$ and $\hat{L}$ of shift-minimal winning and losing coalitions, respectively, to check whether a feasible set of weights exists. If this is not the case we can prune the whole search tree below node $W$. (To check whether $W$ defines a weighted voting game we have to determine the entire set of shift-maximal losing coalitions and solve the corresponding linear program.) We would like to remark that with an increasing number of voters, this simple observation gets increasingly beneficial\footnote{The numbers of weighted voting games and complete simple games coincide for $n\le 5$ voters but their ratio converges to zero with increasing $n$, see also Table~\ref{table_wvg}. An asymptotic upper bound for weighted voting games is given in \cite{keijzer2} and an asymptotic lower bound for complete simple games, there called regular Boolean functions, is given in \cite{0619.05020}.}.

In order to accelerate the check whether a complete simple game is weighted we utilize the following linear program, instead of the one presented in the previous section.
\begin{align}
  &&  \sum_{i=1}^n u_iw_i \ge q &&\forall u\in W,\\
  &&  1+\sum_{i=1}^n v_iw_i \le q &&\forall v\in \widehat{L},\\
  && w_i\ge w_{i-1} && \forall 1\le i\le n-1,\\
  && w_n,q\ge 0.
\end{align}
Thus by using one additional variable $q$, the number of inequalities decreases from $|W|\cdot |\widehat{L}|+n-1$ to $|W|+|\widehat{L}|+n-1$. The key idea to drastically reduce the necessary time for solving the linear programs is, that for a node $W$ and a successor $W'$ the set of variables and many inequalities coincide. So if we use the simplex method, we can perform a warm start using a feasible basis of the linear program corresponding to node $W$ in order to solve the linear program corresponding to node $W'$. Going along this approach drastically reduces the number of performed simplex iterations, i.e.\ base changes, and so the overall running time. 

\begin{table}[htp]
\begin{center}
\begin{tabular}{crrrrrrrrr}
\hline
$\mathbf{n}$&$1$&$2$&$3$&$4$&$5$&$6$&$7$&$8$&$9$\\
$\mathbf{\#}$&$1$&$3$&$8$&$25$&$117$&$1111$&$29373$&$2730164$&$993061482$\\
$\textbf{time}$&&&&&&&$1~s$&$2~m$&$8~d$\\
\hline
\end{tabular}
\caption{Weighted voting games for $n$~voters up to symmetry ($1\sqsupset 2 \sqsupset \dots \sqsupset n$).}
\label{table_wvg}
\end{center}
\end{table}

By combining both ideas and using a re-implementation of the standard simplex method\footnote{\label{fn_da_tautenhahn}We remark that it is also possible to do the enumeration for $9$~voters without the presented ideas, as demonstrated in \cite{da_tautenhahn}. Using some heuristics to find suitable weights on the one hand and to find dual multipliers of the inequalities on the other hand to prove the non-existence of weights, roughly $4$~months of computation time were necessary.}, we were able to enumerate the weighted voting games for $n\le 9$ voters in a reasonable amount of time, see Table~\ref{table_wvg}.

\section{Minimum sum representations}
\label{sec_min_rep}
\noindent
In the following a minimum sum representation of a weighted voting game is a set of integer weights $w_i$ and a quota $q$, such that $\sum\limits_{i=1}^{n} w_i$ is minimal. If additionally $w_i=w_j$ for all $1\le i<j\le n$ with $i~\square~j$, i.e.\ for symmetric voters is required, we speak of a minimum sum representation preserving types, see \cite{root}, \cite{root2} for more details and further definitions of minimum representations. Both minimum representations do exist for each weighted voting game\footnote{Due to the definition of a weighted voting game the corresponding linear program has rational solutions, which can be scaled to be integers.}, but need not be unique if the number $n$ of voters is large enough. So one can ask for a classification of the weighted voting games possessing more than one minimum sum representation (preserving types). In \cite{root} the authors have shown that for at most $7$ voters both minimum sum representations are unique. For $8$~voters the minimum sum representations preserving types are unique and except for $154$~weighted voting games also the minimum sum representations are unique. There are exactly two different minimum sum representations in these $154$~cases.\footnote{We have verified these results using our approach outlined below.}

In order to determine minimum sum representations of weighted voting games and to check whether a given representation is unique we extend our linear program from Section~\ref{sec_partial}. One may simply add the target function $\sum\limits_{i=1}^{n} w_i$ and use integer variables $w_i$. The resulting integer linear programs are that small, that a commercial solver like \texttt{ILOG CPLEX}, or most of the open source alternatives like e.g.\ \texttt{GLPK}, can solve them relatively fast within fractions of a second. But since the number of weighted voting games for $9$~voters is huge (and even more so the number of complete simple games), this is not nearly fast enough. We therefore refrain from using ILPs and use LPs instead. Actually it turned out that the linear program from the previous sections yields integral weights in all cases for at most $7$ voters. For $8$ voters the fractional weights are integral except for $280$~cases, where the denominator equals~$2$.

Suppose that we have integer lower bounds $l_i$ for the weights $w_i$ available (we may use $l_i=0$ at the beginning). In each node $W$ of the search tree, see Section~\ref{sec_partial}, we solve the linear program:

\begin{align}
  \tilde{l}_i:=\min && w_i,\\
  s.t. &&  \sum_{i=1}^n u_iw_i \ge q &&\forall u\in W,\\
  &&  1+\sum_{i=1}^n v_iw_i \le q &&\forall v\in \hat{L},\\
  && w_i\ge w_{i-1} && \forall 1\le i\le n-1,\\
  && w_i\ge l_i && \forall 1\le i\le n,\\
  && w_n,q\ge 0,
\end{align}
for all $i=n,\dots,1$ and update the $l_i$ with the rounded up target values $\left\lceil\tilde{l}_i\right\rceil$ after each solving step, i.e.\ we add some integrality cuts which are valid for every integer representation. Since these $n$ linear programs do not differ much, we highly benefit from warm starting at an optimal basis of the previous solution. In many cases only very few iterations are necessary until the $l_i$ do not change any more. (If $W$ is the entire set of the minimal winning coalitions we have to replace $\hat{L}$ by the set $L$ of all maximal losing coalitions.)

So finally we end up with integral lower bounds $l_i$ for the weights. If the $l_i$ realize the given complete game, then $w_i=l_i$ is a minimum sum representation. Otherwise we store the corresponding complete simple game and the lower bounds $l_i$ as candidates, which have to be treated later on. In the first case the minimum sum representation is unique apart from permutations within equivalence classes of voters.

The big advantage of this approach is that we only have to perform a very small number of simplex iterations, since the consecutive linear programs do not change too much, and that we can pass the lower bounds $l_i$ from a node $W$ to its successor nodes $W'$. Additionally we might use duality of complete simple games, which we have not done.

It turned out that this approach ends up with less than a million candidates out of the $993061482$~weighted voting games for $9$ voters.

Using \texttt{ILOG CPLEX} on the integer linear programming formulation we have determined all minimum sum representations and all minimum 
sum representations preserving types for these candidates separately. There remain exactly $76586$~weighted voting games without a unique 
minimum sum representation\footnote{As a check of the correctness of our computer calculations we have verified that we have found the dual 
games and all examples from the list in \cite{root2}.}. The maximum number of different minimum sum representations for nine voters was three.

\begin{table}[htp]
  \begin{center}
    \begin{tabular}{crrr}
      \hline\\[-2.5mm]
      \textbf{type} & \textbf{2} & \textbf{3} & $\mathbf{\sum}$ \\
      \textbf{1} & 62432 & 624 & 63056\\
      \textbf{2} &     0 & 492 & 492 \\
      \textbf{3} & 12838 &   0 & 12838\\
      \textbf{4} &     0 & 200 & 200\\
      $\mathbf{\sum}$  & 75270 & 1316 & 76586\\
      \hline
    \end{tabular}
    \caption{Number of weighted voting games for $9$ voters without a unique minimum sum representation by type and number of representations.}
    \label{table_counts_1}
  \end{center}
\end{table}

In order to give a more detailed analysis of those examples, see Table~\ref{table_counts_1}, we define several types. As mentioned before 
different minimum sum representations can arise from swaps within an equivalence class of voters, which we call type~1. An example of this kind 
is given by quota $q=49$ and weight vector
\begin{eqnarray*}
  w=\left(24|19|15|8|7{,}7{,}6|2{,}2\right)\!,
\end{eqnarray*}
where we use the $|$'s to abbreviate the equivalence classes $1$, $2$, $3$, $4$, $5~\square~6~\square~7$, and $8~\square~9$. Here we can have 
$3$~different arrangements for the weights $w_5$, $w_6$, and $w_7$.

Instead of a permutation also a redistribution of weights within equivalence classes of voters can occur, which we call type~2. An example 
is given by quota $q=55$ and the weight vectors
\begin{eqnarray*}
  w^{(1)}&=&\left(31|26|23|18|10|7|6|1{,}3\right)\!,\\
  w^{(2)}&=&\left(31|26|23|18|10|7|6|2{,}2 \right)\!,\text{ or}\\
  w^{(3)}&=&\left(31|26|23|18|10|7|6|3{,}1\right)\!.
\end{eqnarray*}

If the weights of the minimum sum representations differ only for voters of different equivalence classes we speak of type~3. These examples 
do not have a unique minimum sum representation preserving types either. By type~4 we denote the remaining cases, where the weights of the 
minimum sum representations differ within and outside of equivalence classes of voters. Such an example is given by the quota $q=56$ and weight 
vectors
\begin{eqnarray*}
  w^{(1)}&=&\left(23|15|13|11|9|8|3|2{,}2\right)\!,\\
  w^{(2)}&=&\left(23|15|13|11|9|8|4|1{,}2\right)\!,\text{ or}\\
  w^{(3)}&=&\left(23|15|13|11|9|8|4|2{,}1\right)\!.
\end{eqnarray*}
In Table~\ref{table_counts_2} we give the number of weighted voting games without a unique minimum sum representation per given number of 
equivalence classes of voters.
\begin{table}[htp]
  \begin{center}
    \begin{tabular}{rrrrrrr}
      \hline
      \textbf{equivalence classes} & 9 & 8 & 7 & 6 & 5 & 4\\
      \textbf{2 representations} & 5718 & 35864 & 24715 & 7659 & 1234 & 80\\
      \textbf{3 representations} &    0 &   402 &   500 &  330 &   76 &  8\\
      \hline
    \end{tabular}
    \caption{Number of weighted voting games for $9$ voters without a unique minimum sum representation by the number of representations and
    equivalence classes of voters.}
    \label{table_counts_2}
  \end{center}
\end{table}

An example for a weighted voting game without a unique minimum sum representation preserving types is given by quota $q=46$ and weights
\begin{eqnarray*}
  w^{(1)}=\left(33|13|12|9|8,8|7|2{,}2\right)\text{ or }
  w^{(2)}=\left(33|13|12|10|8,8|6|2{,}2\right)\!.
\end{eqnarray*}
In total there are $13250$~such weighted voting games, where all of them have exactly two minimum sum representations preserving types. 
The counts by the number of equivalence classes of voters are given in Table~\ref{table_counts_3}. One can easily locate such examples 
for all medium-sized $n\ge 9$. Some additional data is given in Table~\ref{table_additional}.

\begin{table}[htp]
  \begin{center}
    \begin{tabular}{rrrrrr}
      \hline
      \textbf{equivalence classes} & 9 & 8 & 7 & 6 & 5 \\
      \textbf{2 representations} & 5718 & 4992 & 2134 & 392 & 14 \\
      \hline
    \end{tabular}
    \caption{Number of weighted voting games for $9$ voters without a unique minimum sum representation preserving types by the
    number of equivalence classes of voters.}
    \label{table_counts_3}
  \end{center}
\end{table}

see Table~\ref{table_additional}.

\begin{table}[htp]
\begin{center}
\begin{tabular}{lrrrrrrrrr}
\hline
$\mathbf{n}$&$1$&$2$&$3$&$4$&$5$&$6$&$7$&$8$&$9$\\
$\mathbf{\max\min \sum\limits_{i=1}^n w_i}$ &1&2&4&8&15&33&77&202&568\\
$\mathbf{\max\min q}$ &1&2&3&5&9&18&40&105&295\\
$\mathbf{\max\min w_1}$ &1&1&2&3&5&9&18&42&110\\
\hline
\end{tabular}
\caption{Maximum parameters of weighted voting games.}
\label{table_additional}
\end{center}
\end{table}

For $n=9$~voters an example with minimum quota $295$ and minimum sum $568$ is given by the weighted voting game with weights 
$\left(92,84,78,74,67,58,45,40,30\right)$ and $W\!=\!\Big\{\!
110100100,\! 
101011000,\! 
101001011,$ 
$100101101,\! 
100011110,\! 
011110000,\! 
011001101,\! 
010110011,\! 
001111001,\! 
001101110,\! 
000111111\!  
\Big\}$, where we have abbreviated the vector notion of the minimal winning coalitions. An example with $w_1=110$ is given by quota $230$ 
and weights $\left(110,52,48,40,36,28,25,19,7\right)$.

\section{Conclusion}
\noindent
Exhaustively generating all $284\cdot 10^9$ complete simple games for $9$ voters can be done very quickly using a graph theoretic 
representation as cliques in a graph. Solving a linear program in each case to determine whether the game is weighted would be very 
time consuming. A tool to circumvent the exhaustive generation of all complete simple games is the concept of partial complete simple 
games introduced in Section~\ref{sec_partial}. Whether this approach is suitable to enumerate the number of weighted voting games for 
10~voters has to be evaluated in future research. An important application of the enumeration of weighted voting games lies in the inverse 
problem for the power index problem, i.e.\ determine a weighted voting game whose power distribution is close to a given power distribution 
\cite{keijzer}, \cite{keijzer2}, \cite{pre05681536}, \cite{inverse_power_index_problem}.

Concerning the minimum sum representation of weighted voting games we have observed that minimizing the sum of (possibly fractional) 
weights over the inequality system  (5)-(8) very often has integral solutions (at least for a small number of voters). Thus we have 
developed a (simple) cutting plane approach to determine lower bounds on the weights in minimum sum representations. In over 99\% of the 
cases we end up with the unique minimum sum representation. This made it possible to exhaustively generate all minimum sum representations 
(preserving types) of weighted voting games for $9$~voters.

\section*{Acknowledgements}
The author thanks the anonymous referees for carefully reading a preliminary version of this article and giving useful comments to 
improve the presentation.


\end{document}